\input amstex
\frenchspacing
\documentstyle{amsppt}
\magnification=\magstep1
\NoBlackBoxes
\def\FF{\bold F}
\def\Gal{\mathop{\text{\rm Gal}}}
\def\ker{\mathop{\text{\rm ker}}}
\def\ord{\mathop{\text{\rm ord}}}

\def\mod{\bmod}
\def\congr{ \equiv }
\def\O{{\Cal O}}
\def\Que{\bold Q}
\def\Zee{\bold Z}
\def\tilde{\widetilde}
\newcount\refCount
\def\newref#1 {\advance\refCount by 1
\expandafter\edef\csname#1\endcsname{\the\refCount}}
\newref BL   
\newref HAa  
\newref HAb  
\newref HO   
\newref LA   
\newref LG   
\newref LE   
\newref MSa  
\newref MSb  
\newref ST   
\newref STV  
\topmatter
\title 
Prime divisors of the Lagarias sequence
\endtitle
\author Pieter Moree and Peter Stevenhagen\endauthor
\address \hskip-\parindent Mathematisch Instituut,
Universiteit Leiden, P.O. Box 9512, 2300 RA Leiden, The Netherlands
\endaddress
\email moree\@math.leidenuniv.nl,
psh\@math.leidenuniv.nl\endemail

\keywords linear recurrent sequence, Chebotarev density theorem
\endkeywords
\subjclass Primary 11R45; Secondary 11B37
\endsubjclass
\abstract
We solve a 1985 challenge problem posed by Lagarias [\LA]
by determining, under GRH, the density of the set of prime numbers that
occur as divisor of some term of the sequence $\{x_n\}_{n=1}^\infty$ 
defined by the linear recurrence $x_{n+1}=x_n+x_{n-1}$ and the initial
values $x_0=3$ and $x_1=1$.
This is the first example of a `non-torsion' second order recurrent sequence
with irreducible recurrence relation for which we can determine
the associated density of prime divisors.
\endabstract
\endtopmatter

\document

\head 1. Introduction
\endhead

\noindent
In 1985, Lagarias [\LA] showed that the set of prime numbers
that divide some Lucas number has a natural density 2/3 inside
the set of all prime numbers.
Here the Lucas numbers are the terms of the second order recurrent
sequence $\{x_n\}_{n=1}^\infty$ defined by the linear recurrence 
$x_{n+1}=x_n+x_{n-1}$ and the initial values $x_0=2$ and $x_1=1$.
Lagarias's method is a quadratic analogue of the approach
used by Hasse [\HAa, \HAb] in determining, for a given non-zero integer $a$,
the density of the set of the prime divisors of the numbers of the
form $a^n+1$.
Note that the sequence $\{a^n+1\}_{n=1}^\infty$ also satisfies a second order
recurrence.

Hasse and Lagarias apply the Chebotarev density theorem
to a suitable tower of Kummer fields.
Their method of `Chebotarev partitioning'
can be adapted to deal with the class of second order recurrent
sequences that are now known as `torsion sequences' [\BL, \MSa, \STV].
For second order recurrent integer sequences that do not enjoy the rather
special condition of being `torsion', it can no longer be applied.
In the case of the the Lucas numbers, changing the initial values into
$x_0=3$ and $x_1=1$ (while leaving the recurrence $x_{n+1}=x_n+x_{n-1}$
unchanged) leads to a sequence for which Lagarias remarks that his method fails,
and he wonders whether some modification of it can be made to work.

We will explain how non-torsion sequences lead to a question
that is reminiscent of the Artin primitive root conjecture.
In particular, we will see that for a non-torsion sequence, 
there is no number field $F$ (of finite degree) with the property that 
all primes having a given splitting behavior in $F$ divide some
term of the sequence.
It follows that Chebotarev partitioning can not be applied directly.
However, it is possible to combine the technique of Chebotarev partitioning
with the analytic techniques employed by Hooley [\HO] in his proof
(under assumption of the generalized Riemann hypothesis) of 
Artin's primitive root conjecture.
In the case of the modified Lucas sequence proposed by Lagarias,
we give a full analysis of the situation and prove the following
theorem.
\proclaim
{Theorem}
Let $\{x_n\}_{n=0}^\infty$ be the integer sequence defined by
$x_0=3$, $x_1=1$ and the linear recurrence $x_{n+1}=x_n+x_{n-1}$.
If the generalized Riemann hypothesis holds, then the set of prime
numbers that divide some term of this sequence has a natural density.
It equals
$${1573727\over 1569610}\cdot\prod_{p\text{ prime}}
\left(1-{p\over p^3-1}\right)
\approx \ .577470679956\, .
$$
\endproclaim\noindent
Numerically, one finds that 45198 out of the first 78498 primes below $10^6$
divide the sequence: a fraction close to .5758.

\head 2. Second order recurrences
\endhead

\noindent
Let $X=\{x_n\}_{n=0}^\infty$ be a second order recurrent sequence.
It is our aim to determine, whenever it exists,
the density (inside the set of all primes) of the set of
prime numbers $p$ that divide some term of $X$.

We let $x_{n+2}=a_1x_{n+1}+a_0x_n$ be the recurrence satisfied by $X$, and
denote by $f=T^2-a_1T-a_0\in\Zee[T]$ the corresponding 
characteristic polynomial.
We factor~$f$ over an algebraic closure of $\Que$ as
$f=(T-\alpha)(T-\tilde\alpha)$.

In order to avoid trivialities, we will
assume that $X$ does not satisfy a first order recurrence,
so that $\alpha\tilde\alpha=a_0$ does not vanish.
The {\it root quotient\/} $r=r(f)$ of the recurrence, which is only
determined up to inversion, is then defined as $r=\alpha/\tilde\alpha$.
It is either a rational number or a quadratic irrationality of norm 1.
In the {\it separable case\/} $r\ne1$ we have
$$x_n=c\alpha^n+\tilde c\tilde\alpha^n \qquad\hbox{with}\qquad
c={x_1-\tilde\alpha x_0\over \alpha-\tilde\alpha}\quad\hbox{and}\quad
\tilde c={x_1-\alpha x_0\over \tilde\alpha-\alpha}.$$
As our sequence is by assumption not of order smaller than 2,
we have $c\tilde c\ne0$.
Denote by 
$$q={x_1-\alpha x_0\over x_1-\tilde\alpha x_0}=-\tilde c/c
\in\Que(\alpha)^*$$
the {\it initial quotient\/} $q=q(X)$ of $X$.
Just as the root quotient, this is a number determined up to inversion
that is either rational or quadratic of norm 1.
The elementary but fundamental observation for second order
recurrences is that for almost all primes $p$, we have the
fundamental equivalence
$$
p\hbox{ divides }x_n \Longleftrightarrow
-\tilde c/ c=(\alpha/\tilde\alpha)^n \in \O/p\O  \Longleftrightarrow
q=r^n \in (\O/p\O)^*.
$$
Here $\O$ is the ring of integers in the field generated by
the roots of $f$.
This is the ring $\Zee$ if $f$ has rational roots, and the ring of
integers of the quadratic field $\Que[X]/(f)=\Que(\sqrt{a_1^2-4a_0})$
otherwise.
The equivalence above does not make sense for the finitely
many primes $p$ for which either $r$ or $q$ is not invertible
modulo~$p$, but this is irrelevant for density purposes.

In the degenerate case where the root quotient $r$ is a root of unity,
it is easily seen that the set of primes dividing some term of $X$
is either finite or cofinite in the set of all primes.
We will further exclude this case, which includes the
{\it inseparable case\/} $r=1$, for which $q$ is not defined.

As we are essentially interested in the set of
of primes~$p$ for which $q$ is in the subgroup generated by $r$ in the finite
group $(\O/p\O)^*$ of invertible residue classes modulo $p$,
we can formulate the problem we are trying to solve  without
any reference to recurrent sequences.
Depending on whether the root quotient $r$ is rational or quadratic,
this leads to the following.
\proclaim
{Problem 1}
Given two non-zero rational numbers $q$ and $r\ne\pm1$, compute,
whenever it exists, the density of the set of primes $p$ for which we have
$$
 q\mod p\in\langle r\mod p \rangle\subset\FF_p^*.
$$
\endproclaim
\proclaim
{Problem 2}
Let $r$ be a quadratic irrationality of norm $1$ and $\O$ the ring
of integers of $\Que(r)$.
Given an element $q\in\Que(r)$ of norm $1$, compute, whenever it exists,
the density of the set of rational primes $p$ for which we have
$$
 q\mod p\in\langle r\mod p \rangle\subset(\O/p\O)^*.
$$
\endproclaim\noindent
The instances of the two problems above where 
$(q\mod r)$ is a torsion element in the group $\Que(r)^*/\langle r\rangle$
are referred to as torsion cases of the problem, and the
sequences that give rise to them are known as {\it torsion sequences\/}.
The sequences $\{a^n+1\}_{n=0}^\infty$ studied by Hasse, the Lucas sequence
$\{\varepsilon^n+\varepsilon^{-n}\}_{n=0}^\infty$ with
$\varepsilon={1+\sqrt5\over 2}$ treated by Lagarias and the
Lucas-type sequences in [\MSa] are torsion;
in fact, they all have $q=-1$.
The main theorem for torsion sequences, for which we refer to [\STV],
is the following.
\proclaim
{Theorem}
Let $X$ be a second order torsion sequence.
Then the set $\pi_X$ of prime divisors of $X$ has a positive
rational density.
\endproclaim\noindent

\head 3. Non-torsion sequences.
\endhead

\noindent
Problem 1 in the previous section is reminiscent of Artin's famous
question on primitive roots:
given a non-zero rational number $r\ne\pm1$,
for how many primes $p$ does $r$ generate
the group ${\bold F}_p^*$ of units modulo $p$?
(One naturally excludes the finitely many primes $p$ dividing the
numerator or denominator of $r$ from consideration.)
Artin's conjectural answer to this question is based on the observation that
the index $[\FF_p^*:\langle r\rangle]$ is divisible by $j$
if and only if $p$ splits completely in the splitting field
$F_j=\Que(\zeta_j, r^{1/j})$ of the polynomial $X^j-r$ over $\Que$.
Thus, $r$ is a primitive root modulo $p$ if and only $p$ does not split
completely in {\it any\/} of the fields $F_j$ with $j>1$.
For fixed $j$, the set $S_j$ of primes that do split completely in $F_j$
has natural density $1/[F_j:\Que]$ by the Chebotarev density theorem.
Applying an inclusion-exclusion argument to the sets $S_j$,
one expects the set $S=S_1\setminus\cup_{j\ge1} S_j$ of primes for
which $r$ is a primitive root to have natural density
$$
\delta(r)=\sum_{j=1}^\infty {\mu(j)\over [F_j:\Que]}.
\tag{3.1}
$$
Note that the right hand side of (3.1) converges for all 
$r\in\Que^*\setminus \{\pm1\}$ as $[F_j:\Que]$ is 
a divisor of
$\varphi(j)\cdot j$ with cofactor bounded by a constant depending only on $r$.

A `multiplicative version' of the `additive formula' (3.1) for $\delta(r)$
is obtained if one starts from the observation that
$r\in\Que^*\setminus \{\pm1\}$
is a primitive root if and only if $p$ does not split completely
in any field $F_\ell$ with $\ell$ {\it prime\/}. 
The fields $F_\ell$ are of degree $\ell(\ell-1)$ for almost all primes $\ell$,
and using the fact that they are almost `independent', one can
successively eliminate the primes that split completely in some
$F_\ell$ to arrive at a heuristic density
$$
\delta(r)=c_r\cdot\prod_{\ell\text{ prime}} (1-{1\over [F_\ell:\Que]})
=\widetilde c_r\cdot\prod_{\ell\text{ prime}} (1-{1\over \ell(\ell-1)}).
\tag{3.2}
$$
The correction factor $c_r$ for the `dependency' between the fields $F_\ell$
is equal to 1 if the family of
fields $\{F_\ell\}_\ell$ is linearly disjoint over $\Que$, i.\,e.,
if each field $F_{\ell_0}$ is linearly disjoint over $\Que$ from the 
compositum of the fields $F_\ell$ with $\ell\ne\ell_0$.
If $r$ is not a perfect power in $\Que^*$, we have $\widetilde c_r=c_r$.

It turns out that the only possible obstruction to the linear
disjointness of the fields $F_\ell$
occurs when $F_2=\Que(\sqrt r)$ is quadratic of odd discriminant. 
In this case, $F_2$ is contained in the compositum of the fields $F_\ell$
with $\ell$ dividing its discriminant.
The value of $c_r$ is a rational number, and one can derive
a closed formula for it as in [\HO, p. 220].

For example, taking $r=5$, one has $F_2\subset F_5$ and the 
superfluous `Euler factor' $1- [F_5:\Que]^{-1}={19\over 20}$ at $\ell=5$
in the product $\prod_\ell(1- [F_\ell:\Que]^{-1})$ is `removed' by the
correction factor $c_5=\widetilde c_5={20\over 19}$.

It is non-trivial to make the heuristics above into a proof.
As Hooley [\HO] showed, it can be done if one is willing to assume
estimates for the remainder term in the prime number theorem for the
fields $F_j$ that are currently only known to hold under assumption of the
generalized Riemann hypothesis.
One should realize that only when we consider finitely many $\ell$ (or $j$)
at a time, the Chebotarev density theorem gives us the densities we want.
After taking a `limit' over all $\ell$, we only know
that the right hand side of (3.1) or (3.2) is an {\it upper density\/}
for the set of primes $p$ for which $r$ is a primitive root.
We have however no guarantee that we are left with a non-empty set of such $p$.
Put somewhat differently,
we can not obtain primes~$p$ for which $(r\mod p)$ is a primitive root
by imposing a splitting condition on $p$ in a number field $F$ of finite
degree; clearly, there is always some field $F_\ell$ that is linearly
disjoint from $F$, and no splitting condition in $F$ will yield the
`correct' splitting behavior in $F_\ell$.
A similar phenomenon occurs in the analysis of non-torsion cases
of the Problems 1 and 2.
This is exactly what makes non-torsion sequences so much
harder to analyze than torsion sequences.

If $(q\mod r)$ is not a torsion element in $\Que^*/\langle r\rangle$,
then Problem 1 can be treated by a generalization of the arguments used
by Artin.
For each integer $i\ge1$, one considers the set of primes $p$
(not dividing the numerator or denominator of either $q$ or $r$)
for which the index $[\FF_p^*:\langle r\rangle]$ is {\it equal\/} to $i$
and the index $[\FF_p^*:\langle q\rangle]$ is {\it divisible\/}
by $i$.
These are the primes that split completely in the field
$F_{i,1}=\Que(\zeta_i, r^{1/i}, q^{1/i})$, but not in any of the
fields $F_{i,j}=\Que(\zeta_{ij}, r^{1/ij}, q^{1/i})$ with $j>1$.
As before, inclusion-exclusion yields a conjectural value
for the density $\delta_i(r,q)$ of this set of primes, and summing
over $i$ we get
$$
\delta(r,q)=\sum_{i=1}^\infty \delta_i(r,q)\sum_{i=1}^\infty\sum_{j=1}^\infty {\mu(j)\over [F_{i,j}:\Que]}
\tag3.3$$
as a conjectural value for the density in Problem 1.
Note that $\delta_1(r,q)$ is nothing but the primitive root
density $\delta(r)$ from (3.1).

The condition that $(q\mod r)$ is not a torsion element
in $\Que^*/\langle r\rangle$ means
that $q$ and $r$ are multiplicatively independent in $\Que^*$.
In this case $[F_{i,j}:\Que]$ is 
a divisor of $i^2j\cdot\varphi(ij)$ with cofactor bounded by a
constant depending only on $q$ and $r$.
Thus the double sum in (3.3) converges, and under GRH one can
show [\MSb, \ST] that its value is indeed the density one is asked
to determine in Problem 1.

As in Artin's case, one can obtain a multiplicative version of (3.3)
by a `prime-wise' approach.
One notes that the inclusion of subgroups
$\langle q\mod p\rangle\subset\langle r\mod p \rangle$ in $\FF_p^*$
means that for all primes $\ell$, we have an inclusion
$$
\langle q\mod p\rangle_\ell\subset\langle r\mod p \rangle_\ell
\tag3.4
$$
of the $\ell$-primary parts of these subgroups.
If we fix both $\ell$ and the number $k=\ord_\ell(p-1)$ of factors $\ell$
in the order of $\FF_p^*$, this condition can be rephrased in terms of the
splitting behavior of $p$ in the number field
$$
\Omega_\ell^{(k)}=\Que(\zeta_{\ell^{k+1}}, r^{1/\ell^k}, q^{1/\ell^k}).
\tag{3.5}
$$
More precisely, we have $\ord_\ell(p-1)=k$ if and only if
$p$ splits completely in $\Que(\zeta_{\ell^k})$ but not in
$\Que(\zeta_{\ell^{k+1}})$;
of the primes $p$ that meet this condition, we want those $p$ for which
the order of the Frobenius elements over $p$ in
$\Gal(\Que(\zeta_{\ell^k},q^{1/\ell^k})/\Que(\zeta_{\ell^k}))$ divides
the order of the Frobenius elements over $p$ in
$\Gal(\Que(\zeta_{\ell^k},r^{1/\ell^k})/\Que(\zeta_{\ell^k}))$.
By the Chebotarev density theorem, one finds that the set of primes $p$
with $\ord_\ell(p-1)=k$, which has density $\ell^{-k}$ for $k\ge1$,
is a union of two sets that each have a density:
the set of primes $p$ for which the inclusion (3.4) holds
and the set of $p$ for which it does not.
This `Chebotarev partitioning' allows us to compute, for each $\ell$,
the density of the primes $p$ for which we have the inclusion (3.4):
summing over $k$ in the previous argument yields a lower density, and this
is the required density as we can apply the same argument to the
complementary set of primes.

For all but finitely many $\ell$, the fields
$\Que(\zeta_{\ell^{k+1}},q^{1/\ell^k})$ and
$\Que(\zeta_{\ell^{k+1}},r^{1/\ell^k})$
are linearly disjoint extensions of $\Que(\zeta_{\ell^{k+1}})$ with
Galois group $\Zee/\ell^k\Zee$ for all $k\ge0$.
In this case the set of primes $p$ with $\ord_\ell(p-1)=k$ violating (3.4)
has density 
$$\ell^{-k} \sum_{i=1}^k \ell^{-i} (\ell^{-(i-1)} - \ell^{-i})
=(\ell^{-k}-\ell^{-3k})/(\ell+1).$$
Summing over $k$, we
find that (3.4) does not hold for a set of primes of density
$\ell/(\ell^3-1)$.
As the fields
$$
\Omega_\ell=\textstyle\bigcup_k \Omega_\ell^{(k)}\Que(\zeta_{\ell^\infty},r^{1/\ell^\infty}, q^{1/\ell^\infty})
$$
for prime values of $\ell$ form a linearly disjoint family if we exclude
finitely many `bad' primes $\ell$, the multiplicative analogue of (3.3)
reads
$$
\delta(r, q)c_{q,r}\cdot\prod_{\ell\text{ prime}} \left(1-{\ell\over \ell^3-1}\right).
\tag{3.6}
$$
As is shown in [\MSb], the `correction factor' $c_{q,r}$ is 
a rational number that admits a somewhat involved description
in terms of $q$ and $r$.
In practice, one finds its value most easily by starting from
the additive formula (3.3).
\medskip\noindent
In the situation of Problem 2, the arguments just given
can be taken over without 
substantial changes from the rational case
when one restricts to those rational primes
$p$ that split completely in $\O$.
Writing $K=\Que(r)$ and
$$F_{i, j}=K(\zeta_{ij}, r^{1/ij}, q^{1/i}),$$
we find that the density of the rational primes $p$ that are split in $\O$
and for which we have $ q\mod p\in\langle r\mod p \rangle\subset(\O/p\O)^*$
equals
$$\delta_{\text{split}}(r, q)\sum_{i=1}^\infty\sum_{j=1}^\infty {\mu(j)\over [F_{i,j}:\Que]};
\tag{3.7}$$
as in the case of Problem 1, one needs to assume the validity
of the generalized Riemann hypothesis for this result.
By (3.7), the computation of $\delta_{\text{split}}$ amounts
to a degree computation for the family of fields $\{F_{i,j}\}_{i,j}$.
In fact, because of the numerator $\mu(j)$ in (3.7)
one may restrict to the case where $j$ is squarefree.
As in the case of Problem 1, one finds (under GRH) that the split density equals
$$\delta_{\text{split}}(r, q)c_{q,r}^+\cdot\prod_{\ell\text{ prime}} \left(1-{\ell\over \ell^3-1}\right)
\tag{3.8}$$
for some rational number $c_{q,r}^+$.
The next section provides a typical example of such a computation.
It shows that the value of $c_{q,r}^+$ is not as simple a fraction
as the analogous factor $c_r$ in (3.2).

For the rational primes $p$ that are inert in $\O$, the determination
of the corresponding density $\delta_{\text{inert}}$ is more involved
than in the split case.
The group $(\O/p\O)^*$ in Problem 2 is now cyclic of order $p^2-1$,
and $(q\mod p)$ and $(r\mod p)$ are elements of the kernel
$$
\kappa_p=\ker[N: (\O/p\O)^*\longrightarrow \FF_p]
\tag{3.9}
$$
of the norm map, which is cyclic of order $p+1$.
In order to have the inclusion (3.4) of subgroups of $\kappa_p$ for
all primes $\ell$,
we fix $\ell$ and $k=\ord_\ell(p+1)\ge1$ and rephrase (3.4) in terms of
the splitting behavior of $p$ in the quadratic counterpart
$$
\Omega_\ell^{(k)}=K(\zeta_{\ell^{k+1}}, r^{1/\ell^k}, q^{1/\ell^k})
$$
of (3.5).
Let us assume for simplicity that $\ell$ is an odd prime,
and that $K$ is not the quadratic subfield of $\Que(\zeta_\ell)$.
Then the requirement that $p$ be inert in $K$ and satisfy
$\ord_\ell(p+1)=k\ge1$ means that the Frobenius element of $p$ in
$\Gal(K(\zeta_{\ell^{k+1}})/\Que)$ is non-trivial on $K$ and has order~$2\ell$ 
when restricted to $\Que(\zeta_{\ell^{k+1}})$.
Let $B_k\subset K(\zeta_{\ell^{k+1}})$ be the fixed field
of the subgroup generated by such a Frobenius element.
Then $B_k$ does not contain $K$ or $\Que(\zeta_\ell)$, and
$B_k\subset K(\zeta_{\ell^k})$ is a quadratic extension.
Let $\sigma_k$ be the non-trivial automorphism of this extension.
Then $\sigma_k$ acts by inversion on $\zeta_{\ell^k}$, 
and the norm-1-condition on $q$ and $r$ 
means that $\sigma_k$ also acts by inversion on~$q$ and~$r$.
The Galois equivariancy of the Kummer pairing 
$$
\Gal(K(\zeta_{\ell^k},r^{1/\ell^k}, q^{1/\ell^k})/K(\zeta_{\ell^k}))
\times \langle q, r\rangle\longrightarrow
\langle \zeta_{\ell^k}\rangle$$
shows that the natural action of $\sigma_k$ on
$\Gal(K(\zeta_{\ell^k},r^{1/\ell^k}, q^{1/\ell^k})/K(\zeta_{\ell^k}))$
is trivial, so $\Omega_\ell^{(k)}$ is abelian over $B_k$.
It is the linearly disjoint compositum of the cyclotomic extension
$B_k\subset B_k(\zeta_{\ell^{k+1}})$ and the abelian extension 
$$B_k\subset
B_k(r^{1/\ell^k}+r^{-1/\ell^k},q^{1/\ell^k}+q^{-1/\ell^k}).$$
For almost all $\ell$, the group $\Gal(\Omega_\ell^{(k)}/B_k)$
is isomorphic to $\Zee/2\ell\Zee\times (\Zee/\ell^k\/\Zee)^2$.

Just like in the rational case, we want those primes $p$ that have splitting 
field $B_k$ inside $K(\zeta_{\ell^{k+1}})$ and for which the order
of the Frobenius elements over $p$ in
$\Gal(B_k(q^{1/\ell^k}+q^{-1/\ell^k})/B_k)$ divides
the order of the Frobenius elements over $p$ in
$\Gal(B_k(r^{1/\ell^k}+r^{-1/\ell^k})/B_k)$.
By the Chebotarev partition argument, we find again that for a
`generic' prime $\ell$, a fraction $\ell/(\ell^3-1)$ of the 
primes $p$ that are inert in $\O$ violates (3.4).
Here `generic' means that $\ell$ is odd and that $\Omega_\ell^{(k)}$ has degree 
$2\ell^{3k}(\ell-1)$ for $k\ge1$.
Under GRH, one can again deduce that the inert density
$\delta_{\text{inert}}(q, r)$ equals a rational constant
$c_{q, r}^-$ times the infinite Euler product occurring in (3.6) and (3.8).

In general, there are various subtleties that need to be taken care of in the
analysis above for $\ell=2$, when 2-power roots of $q$ and $r$ are 
adjoined to $K$ or $K(\zeta_{2^n})$. 
We do not go into them in this paper.
In the example in the next section, we deal with these complications by 
combining a simple ad hoc argument for a few 
`bad' $\ell$ with the standard treatment for the `good' $\ell$.

\head 4. The Lagarias example
\endhead

\noindent
We now treat the explicit example of the modified Lucas sequence which
is the subject of the theorem stated in the introduction.
The roots of the characteristic polynomial $X^2-X-1$ of the recurrence 
are $\varepsilon={1+\sqrt 5\over 2}$ and its conjugate
$\tilde\varepsilon={1-\sqrt 5\over 2}$.
The initial values $x_0=3$ and $x_1=1$ yield an initial quotient
$q={1-3\varepsilon\over 1-3\tilde\varepsilon}$ of the sequence.
As $\pi_{11}=1-3\varepsilon\in \O=\Zee[\varepsilon]$ has norm $-11$,
we find that we have to solve Problem 2 for
$$
q={\pi_{11}\over\tilde\pi_{11}}={\pi_{11}^2\over-11}
\qquad\text{and}\qquad
r={\varepsilon\over\tilde\varepsilon}=-\varepsilon^2.
$$
We set $K=\Que(\varepsilon)=\Que(\sqrt 5)$ and
$F_{i, j}=K(\zeta_{ij}, r^{1/ij}, q^{1/i})$
as in the previous section.
\proclaim
{4.1. Lemma}
For $i, j\in\Zee_{\ge1}$ we have
$[F_{i,j}:\Que]=2^{1-t} i^2 j\varphi(ij)$,
with
$$t=t_{i,j}=\cases
\#\{d\in \{4, 5, 11\}: d|ij\}&       \text{if $i$ is even;}\\
\#\{d\in \{4, 5\}: d|ij\}&       \text{if $i$ is odd and $j$ is even;}\\
\#\{d\in \{5\}: d|ij\}&       \text{if $ij$ is odd.} \endcases
$$
\endproclaim\noindent
{\bf Proof.\/}
As $K=\Que(\sqrt 5)$ is the quadratic subfield of $\Que(\zeta_5)$,
the field $K(\zeta_{ij})$ has degree $2\varphi(ij)$ over $\Que$
if 5 does not divide $ij$ and degree $\varphi(ij)$ if it does.

As $q={\pi_{11}\over\tilde\pi_{11}}$ is a quotient of two
non-associate prime elements in $\O$ and $\varepsilon$ is a fundamental unit
in $\O$, the polynomials $X^i-q$ and $X^{ij}-r$ are irreducible in $K[X]$
for all $i, j\in \Zee_{\ge1}$ by a standard result as [\LG, Theorem VI.9.1].
Moreover, the extension $K\subset K(q^{1/i})$ generated by a zero of $X^i-q$
is totally ramified at the primes of $K$ lying over 11, whereas the
extension $K\subset K(r^{1/ij})$ generated by a zero of $X^{ij}-r$
is unramified above 11.
It follows that $K\subset K(q^{1/i}, r^{1/ij})$ is of degree $i^2j$ 
for all $i, j\in \Zee_{\ge1}$.

The intersection $K(q^{1/i}, r^{1/ij})\cap K(\zeta_{ij})$ is contained
in the maximal abelian subfield $K_0$ of $K(q^{1/i}, r^{1/ij})$,
which equals 
$$
K_0=\cases
K(\sqrt q, \sqrt r)=K(\sqrt{-11}, \zeta_4)&       \text{if $i$ is even;}\\
K(\sqrt r)=K(\zeta_4)&       \text{if $i$ is odd and $j$ is even;}\\
K&       \text{if $ij$ is odd.} \endcases
$$
One trivially computes
$K_0\cap K(\zeta_{ij})$, and the lemma follows.\hfill$\square$
\medskip\noindent
We will need the preceding lemma only for squarefree $j$.
In this case, we simply have $t=\#\{d\in \{5\}: d|ij\}$ for odd $i$.

If we substitute the explicit degrees from Lemma 4.1 in (3.7), we find
that the split density for our example equals
$$
\delta_{\text{split}}{1\over 2} \sum_{i=1}^\infty\sum_{j=1}^\infty
2^{t_{i,j}}{\mu(j)\over i^2j\varphi(ij)},
$$
where $t_{i,j}$ is as in Lemma 4.1.
If we set
$$
S_{m, n}=\sum_{{i=1\atop m|i}}^\infty \sum_{{j=1\atop mn|ij}}^\infty
{\mu(j)\over i^2j\varphi(ij)},$$
then the expression above may be rewritten as
$$
2\delta_{\text{split}}S_{1,1}+S_{2,2}+S_{2,5}+S_{2,11}+S_{2,10}+S_{2,22}+S_{2,55}+S_{2,110}
+(S_{1,5}-S_{2,5}).
$$
It is elementary to show [\MSb, Theorem 4.2] that $S_{m, n}$ is the
rational multiple
$$S_{m, n}{S\over m^3n^3} \prod_{p|n} {-p^4\over p^3-p-1}
                \prod_{{p|m\atop p\nmid n}} {p^3+p^2\over p^3-p-1}$$
of the universal constant $S=\prod_{p\text{ prime}} (1-{p\over p^3-1})$
for $m, n\in\Zee_{\ge1}$.
Simple arithmetic now yields the value
$$\delta_{\text{split}}{712671\over 1569610}\cdot S
\approx \ .26151
\tag{4.2}$$
for the density (under GRH) of the primes $p\congr\pm1 \mod 5$
dividing the Lagarias sequence.
Numerically, one finds that 20416 primes out of the 78498 primes
below $10^6$ are split in $K$ and divide our sequence: a fraction
close to $.2601$.

For the inert primes of $K$, which satisfy $p\congr\pm2\mod 5$,
we have a closer look at the
`bad' primes 2, 5 and~11.
As for the rational problem, we define the extensions
$$
\Omega_\ell=K(\zeta_{\ell^\infty},r^{1/\ell^\infty}, q^{1/\ell^\infty})
$$
for primes $\ell$ and note that, by Lemma 4.1, the family consisting of
the extensions $\Omega_2\Omega_5\Omega_{11}$ and
$\{\Omega_\ell\}_{\ell\ne 2, 5, 11}$ of $K$ is linearly independent over $K$.

Our first observation is that for the inert primes $p$, the order
$p+1$ of the group $\kappa_p$ in (3.9) is never divisible by 5.
Condition (3.4) is therefore automatic for the prime $\ell=5$, and we can
disregard the splitting behavior of $p$ in $\Omega_5$.

We next observe that for inert $p$, the element $r=-\varepsilon^2$ satisfies
$$
r^{(p+1)/2}\congr(-1)^{(p-1)/2}\mod p.
\tag{4.3}$$
For primes $p\congr 3\mod 4$, this shows that 
$\langle r\mod p\rangle_2$ is the 2-Sylow subgroup of $\kappa_p$,
so that (3.4) is again automatic for $\ell=2$.
When we now impose that the inert primes congruent to $3\mod 4$,
which form a set of primes of density $1/4$, have the correct
splitting behavior in the extensions $\Omega_\ell$ for $\ell\ne 2,5$,
we are dealing with a linearly disjoint family and find (under GRH) 
that the set of these primes has density
$$
{1\over 4} \cdot\prod_{\ell\ne 2, 5} \left(1-{\ell\over \ell^3-1}\right)
={1\over 4} \cdot {7\over 5}\cdot{124\over 119} \cdot S.
$$
We next consider the inert primes $p\congr1\mod 4$.
For these $p$, the congruence (4.3)
shows that $(r\mod p)$ has odd order in $\kappa_p$, so (3.4)
is satisfied for $\ell=2$ if and only if 
the order of $q=-\pi_{11}^2/11$ in $\kappa_p$ is also odd.
As $\kappa$ is a cyclic group of order $p+1\congr2\mod 4$,
the order of $\bar q=(q\mod p)$ is odd if and only if $\bar q$
is a square in~$\kappa_p$.
Let $x\in \O/p\O$ be a square root of $\bar q$.
If $x$ is in $\kappa_p$, i.e., if $x$ has norm~1 in~$\FF_p$,
then its trace $x+1/x$ is in $\FF_p$, and we find that
$$(x+{1\over x})^2=2+\bar q+{1\over \bar q}={1\over -11}\mod p$$
is a square modulo $p$.
If $x$ is not in $\kappa_p$, then $x$ has norm $-1$ in $\FF_p$ and
$$(x-{1\over x})^2=2-\bar q-{1\over \bar q}={3^2\cdot 5\over -11}\mod p$$
is a square modulo $p$.
As 5 is not a square modulo our inert prime $p$, we deduce
$$
\text{$(q\mod p)$ has odd order in $\kappa_p$}
\quad \Longleftrightarrow\quad
\text{$-11$ is a square modulo $p$}.
\tag{4.4}
$$
\noindent
If $p$ satisfies the equivalent conditions of (4.4),
then $p$ is a square modulo 11
by quadratic reciprocity, and we have $11\nmid p+1$.
It follows that in this case, (3.4) is satisfied for $\ell=2, 5$ and $11$.
Thus, the set of inert primes $p\congr 1\mod 4$ satisfying the
quadratic condition (4.4) is a set of primes of density 1/8, and
the subset of those $p$ that have the correct
splitting behavior in the extensions $\Omega_\ell$ for $\ell\ne 2, 5, 11$
has (under GRH) density
$$
{1\over 8} \cdot\prod_{\ell\ne 2, 5, 11} \left(1-{\ell\over \ell^3-1}\right)
={1\over 4} \cdot {7\over 5}\cdot{124\over 119} \cdot {1330\over 1319}\cdot S.
$$
Adding the fractions obtained for the inert primes congruent to
$3\mod 4$ and to $1\mod 4$, we obtain
$$
\delta_{\text{inert}}=
{61504\over 112115}\cdot S \approx 0.3159598798268.
$$
Numerically, one finds that 24781 primes out of the 78498 primes
below $10^6$ are inert in $K$ and divide our sequence: a fraction
close to $.3157$.

The sum $\delta_{\text{split}}+\delta_{\text{inert}}$ is the
value
$({712671\over 1569610} +{61504\over 112115})\cdot S{1573727\over 1569610}\cdot S$
mentioned in the theorem in the introduction.

\enddocument